\documentclass[11pt,twoside]{article}
\usepackage{amsmath}
\usepackage{amsthm}
\usepackage{amssymb}
\usepackage{array}
\usepackage{geometry}
\usepackage{graphicx}
\usepackage{enumerate}
\usepackage{lscape,graphicx}
\usepackage{fancyhdr}
\usepackage{makeidx}
\usepackage{amsfonts}
\usepackage{graphicx}
\usepackage{CJK}
\usepackage{amsmath}
\usepackage{latexsym}
\usepackage{amssymb}
\usepackage{mathrsfs}
\usepackage{amsthm}
\usepackage{array}
\usepackage{enumerate}
\usepackage{longtable}
\usepackage{multirow}
\makeindex
\usepackage{wasysym}
\usepackage{chapterbib}

\normalsize
\begin{document}


\newpage
\setcounter{equation}{0}
\begin{center}
\vskip1cm
{\Large
\textbf{Different Estimation Procedures For Topp Leone Exponential And Topp Leone q Exponential Distruibution} }
\vskip.5cm
%

\vspace{0.50cm}
\textbf{\large Nicy Sebastian$^{1}$  and Rajitha V R$^{2}$ }\\
$^{1}$Department of Statistics,
St.Thomas College, Thrissur, India -680 001\\
{Email:{\tt{nicycms@gmail.com}}}\\
$^{2}$Department of Statistics,
St. Aloysius College
Elthuruth, Thrissur, 
 India- 680 611\\ {Email:{\tt{rajithavr95@gmail.com}}}\\
\end{center}

\thispagestyle{empty}
\begin{center} {\bf }  \vskip 0.10truecm \noindent \\

\vskip.5cm\noindent \centerline{\bf Abstract}
\end{center}
\par
Topp Leone q Exponential Distruibution is a continuous model distribution used for modelling lifetime phenomena. In this study, we introduce different estimation methods for the unknown parameters of Topp Leone Exponential(TLE) distribution and Topp Leone q Exponential(TLqE) distribution. \\



\vskip.5cm{\section*{1.\hskip.3cm Introduction}}
\vskip.3cm
In the survival, a number of continuous univariate distributions have been widely
utilized for demonstrating information in numerous areas, for example, biology,
medicine, engineering, public health, epidemiology and economics. In any case, applied areas, for example, lifetime analysis obviously require expanded types of these
distributions. In this way, a few classes of distributions have been built by extending common families of continuous distributions.
These generalized distributions
give greater adaptability by including \textquoteleft at least one\textquoteright parameters to the standard
model. The Topp-Leone(TL) distribution was introduced by Topp and Leone in
1955 (Topp and Leone, 1955). Topp-Leone Generalized(TLG) family of distributions was inferred by Rezaei et al. (2016). The distribution and density function of
proposed family is known by
\begin{eqnarray}
F_{TLG}(x)&=& 2\alpha \int_{0}^{G(X)} t^{\alpha -1}(1-t)(2-t)^{\alpha-1}dt \nonumber \\ &=& G(x)^\alpha (2-G(x))^\alpha
\end{eqnarray}
differentiating, we get the corresponding pdf,
\begin{equation}
f_{TLG}(x)=2\alpha g(x)(1-G(x))G((x)^{\alpha-1}(2-G(x))^{\alpha-1})
\end{equation}
where $G(x) $ and $ g(x)$ are cdf and pdf respectively ,and $\alpha\geq 0$.
\vskip.3cm {\subsection*{1.1.\hskip.3cm Topp Leone Exponential distribution(TLE)}}

\paragraph{}
In reliability analysis, a frequently used distribution is exponential distribution(Crowder et. al., 1994). Its characterizing property is its constant hazard
function. Due to this, exponential distribution is sometimes not suitable for
analyzing data. This implies the need for more generalization. In such situations we use distribution called Topp-Leone Exponential distribution (TLE).
TLE distribution comes as the combination of TL distribution and exponential distribution. Here TL distribution is the generator and exponential is
the parent distribution
\\ For creating the TLE, we need cdf $G(x)$ and pdf $g(x)$ of exponential distribution
\begin{equation}
G(x)=1-\exp(-\lambda x);x\geq 0,\lambda \geq 0
\end{equation}
and
\begin{equation}
g(x)=\lambda \exp (-\lambda x)
\end{equation}

The TLE distribution is obtained by taking  and  into  and .
\\ A random variable X possessing TLE distribution has the cdf and the pdf
\begin{eqnarray}
F_{TLE} (x)&=& (1-\exp (-\lambda x)
)^\alpha (2-(1-\exp (-\lambda x)))^\alpha
\\ &=& (1-\exp(-2\lambda x))^\alpha
\end{eqnarray}
and
\begin{equation}
f_{TLE} (x)=2 \alpha \lambda \exp (-2\lambda x)(1-\exp (-2 \lambda x))^{\alpha-1}
\end{equation}
here $\alpha$ is the shape parameter and $\lambda $ is the scale parameter.
\vskip.3cm {\subsection*{1.2.\hskip.3cm Topp Leone $q$ Exponential distribution(TL$q$E)}}
\paragraph{}
$q$-exponential distribution is a higher version of exponential distribution. It
provides more flexibility regarding to its decay than exponential distribution.
It does not have the limitation of constant hazard rate. And thus allows modeling system improvement and degeneration
\\ Here we generalize q-exponential distribution into Topp Leone q-Exponential
(TLqE) distribution by exponential distribution. The TL distribution combined with q-exponential distribution gives TLqE distribution. Here q-exponential
distribution is the parent distribution and TL distribution is the generator
distribution. The distribution function and density function of q-exponential
distribution is given as
\begin{equation}
F_{qE} (X)=1-[1-(1-q)\lambda x ]^{\frac{2-q}{1-q}}
\end{equation}
and
\begin{equation}
f_{qE} (x)=(2-q)\lambda [1-(1-q)\lambda x]^{\frac{1}{1-q}} ;x\geq 0, \lambda \geq 0 ,q\leq 2,q\neq 0
\end{equation}
respectively.Substituting these in the cdf and pdf of TLG distribution we
get the cdf and pdf of TLqE distribution.

A random variable X has a TLqE distribution, if X has the cdf and pdf
respectively as follows,
\begin{eqnarray}
F_{TLqE}(x)&=&[1-[1-(1-q)\lambda x]^{\frac{2-q}{1-q}}]^\alpha [2-[1-[1-(1-q)\lambda x ]^{\frac{2-q}{1-q}}]]^\alpha
\\ &=&[1-[1-(1-q)\lambda x ]^{2{\frac{2-q}{1-q}}}]^\alpha
\end{eqnarray}
$x\geq 0,\lambda , \alpha \geq 0,q \leq 2 ,q\neq 0 $
\\ and
\begin{eqnarray}
f_{TqE}(x)&=&2\alpha (2-q)\lambda[1-(1-q)\lambda x]^{\frac{1}{1-q}} [1-(1-q)\lambda x] ^{2 {\frac{2-q}{1-q}}}
\nonumber	\\ &&[1-[1-(1-q)\lambda x]^{2 {\frac{2-q}{1-q}}}]^{\alpha -1} \nonumber
\\ &=& 2\alpha \lambda (2-q) [1-(1-q)\lambda x ]^{\frac{3-q}{1-q}} [1-[1-(1-q)\lambda  x]^{2{\frac{2-q}{1-q}}}]^{\alpha -1}
\end{eqnarray}
\vskip.3cm {\subsection*{2.\hskip.3cm Parameter Estimation}}
	\paragraph{}We now explore the statistical aspect of the TLE distribution and TLqE distribution and investigate the estimation of the unknown parameters $\lambda$,$\alpha$ and $q$ by five methods. Howafter $x_1,x_2,...x_n$ denote realizations from a random sample of size $n$ from X and $x_{(1)},x_{(2)},...,x_{(n)}$ their ascending order.
	\vskip.3cm {\subsection*{2.1.\hskip.3cm Topp Leone Exponential Distribution}}
	\vskip.3cm {\subsection*{2.1.1\hskip.3cm Method of Least Squares Estimation}}
	\paragraph{}Here we consider the least squares estimation introduced by in their joint work. The least square estimates(LSEs) of $\lambda$ and $\alpha$ can be determined by minimizing the least square function, with respect to $\lambda$ and $\alpha$. The least square function is defined by
	\begin{eqnarray}
	LS(\lambda,\alpha)&=&\sum_{i=1}^{n}\big[F(x_{(i)};\lambda,\alpha)-\dfrac{i}{n+1}\big]^2  \\ &=& \sum_{i=1}^{n}\big[ [1-e^{-2\lambda x_{(i)}}]^\alpha-\dfrac{i}{n+1}\big]^2 \nonumber
	\end{eqnarray}
	Thus, the least square estimates can be obtained by solving the following equations simultaneously: $\partial LS(\lambda,\alpha)/\partial \lambda =0$ and $\partial LS(\lambda,\alpha)/\partial \alpha =0$, where
	\begin{eqnarray}
	\dfrac{\partial LS(\lambda,\alpha)}{\partial \alpha}&=&2\sum_{i=1}^{n}\big[ [1-e^{-2\lambda x_{(i)}}]^\alpha-\dfrac{i}{n+1}\big][1-e^{-2\lambda x_{(i)}}]^\alpha \log{[1-e^{-2\lambda x_{(i)}}]} \nonumber \\
	\end{eqnarray}
	and
	\begin{eqnarray}
	\dfrac{\partial LS(\lambda,\alpha)}{\partial \lambda}&=&2\sum_{i=1}^{n}\big[ [1-e^{-2\lambda x_{(i)}}]^\alpha-\dfrac{i}{n+1}\big] 2\alpha x_{(i)} [1-e^{-2\lambda x_{(i)}}]^{\alpha-1}{e^{-2\lambda x_{(i)}}} \nonumber \\
	\end{eqnarray}
	\vskip.3cm {\subsection*{2.1.2.\hskip.3cm Method of Weighted Least Squares estimation}}
	\paragraph{}The weighted least square estimates(WLSEs) of $\lambda$ and $\alpha$ can be obtained by minimizing the weighted least square function, with respect to $\lambda$ and $\alpha$. The weighted least square function is defined by
	\begin{eqnarray}
	LSW(\lambda,\alpha)&=&\sum_{i=1}^{n}\dfrac{(n+1)^2(n+2)}{i(n-i+1)}\big[F(x_{(i)};\lambda,\alpha)-\dfrac{i}{n+1}\big]^2  \\ &=& \sum_{i=1}^{n}\big[ [1-e^{-2\lambda x_{(i)}}]^\alpha-\dfrac{i}{n+1}\big]^2 \nonumber
	\end{eqnarray}
	Thus, the weighted least square estimates can be obtained by solving the following equations simultaneously: $\partial LSW(\lambda,\alpha)/\partial \lambda =0$ and $\dfrac{\partial LSW(\lambda,\alpha)}{\partial \alpha}=0,$ where
	\begin{eqnarray}
	&&\dfrac{\partial LSW(\lambda,\alpha)}{\partial \alpha}\nonumber\\
&&\nonumber\\
& &\quad =2\sum_{i=1}^{n}\dfrac{(n+1)^2(n+2)}{i(n-i+1)}\big[ [1-e^{-2\lambda x_{(i)}}]^\alpha-\dfrac{i}{n+1}\big][1-e^{-2\lambda x_{(i)}}]^\alpha \log{[1-e^{-2\lambda x_{(i)}}]} \nonumber\\
	\end{eqnarray}
	and
	\begin{eqnarray}
	&&\dfrac{\partial LSW(\lambda,\alpha)}{\partial \lambda}\nonumber\\
&&\nonumber\\
& &\quad =2\sum_{i=1}^{n}\dfrac{(n+1)^2(n+2)}{i(n-i+1)}\big[ [1-e^{-2\lambda x_{(i)}}]^\alpha-\dfrac{i}{n+1}\big]  2\alpha x_{(i)} [1-e^{-2\lambda x_{(i)}}]^{\alpha-1}{e^{-2\lambda x_{(i)}}}\nonumber\\
	\end{eqnarray}
	\vskip.3cm {\subsection*{2.1.3\hskip.3cm Method of Cramer-von Mises Distance Estimation}}
	\paragraph{}The Cramer-von Mises estimator (CME) is a type of minimum
	distance estimators (also called maximum goodness of
	fit estimators) which is based on the difference between
	the estimate of the cumulative distribution function and the
	empirical distribution function .
	MacDonald  motivates the choice of Cramer-von
	Mises type minimum distance estimators providing empirical
	evidence that the bias of the estimator is smaller than the
	other minimum distance estimators.The cramer-von Mises minimum distance estimates(CVEs) of $\alpha$ and $\lambda$ is determined by minimizing the cramer-von Mises minimum distance function , with respect to $\lambda$ and $alpha$. The cramer-von Mises minimum distance function is defined by
	\begin{eqnarray}
	C(\lambda,\alpha)&=&\frac{1}{12n}+\sum_{i=1}^{n}\big[F(x_{(i)};\lambda,\alpha)-\dfrac{2i-1}{2n}\big]^2 \\ &=&\frac{1}{12n}+\sum_{i=1}^{n}\big[[1-e^{-2\lambda x_{(i)}}]^\alpha-\dfrac{2i-1}{2n}\big]^2 \nonumber
	\end{eqnarray}
	
	Thus, the cramer-von Mises minimum distance estimates can be obtained by solving the following equations simultaneously: $\partial C(\lambda,\alpha)/\partial \lambda =0$ and $\partial C(\lambda,\alpha)/\partial \alpha =0$, where
	\begin{eqnarray}
	\dfrac{\partial C(\lambda,\alpha)}{\partial \alpha}&=&\sum_{i=1}^{n}\big[[1-e^{-2\lambda x_{(i)}}]^\alpha-\dfrac{2i-1}{2n}\big][1-e^{-2\lambda x_{(i)}}]^\alpha \log{[1-e^{-2\lambda x_{(i)}}]} \nonumber \\
	\end{eqnarray}
	and
	\begin{eqnarray}
	\dfrac{\partial C(\lambda,\alpha)}{\partial \lambda}&=&\sum_{i=1}^{n}\big[[1-e^{-2\lambda x_{(i)}}]^\alpha-\dfrac{2i-1}{2n}\big]2 \alpha x_{(i)}[1-e^{-2\lambda x_{(i)}}]^{\alpha-1}e^{-2\lambda x_{(i)}} \nonumber \\
	\end{eqnarray}
\vskip.3cm {\subsection*{2.1.4\hskip.3cm Method of Anderson-Darling Estimation}
	\paragraph{} The method of Anderson-Darling estimation was introduced by Anderson-Darling in the context of statistical tests. In TLE distribution, the Anderson-Darling estimates(ADEs) of$\lambda$ and $\alpha$ can be determined by minimizing the Anderson-Darling function, with respect to $\lambda$ and $\alpha$. The Anderson-Darling function is defined by
	\begin{eqnarray}
	A(\lambda,\alpha)&=&-n-\frac{1}{n}\sum_{i=1}^{n}(2i-1)\{\log[F(x_{(i)};\lambda,\alpha)]+\log[1-F(x_{(n+1-i)};\lambda,\alpha)]\} \nonumber \\ \\ &=&-n-\frac{1}{n}\sum_{i=1}^{n}(2i-1)\{\log[1-e^{-2\lambda x_{(i)}}]^\alpha+\log[1-e^{-2\lambda x_{(n+1-i)}}]^\alpha\}
	\end{eqnarray}
	Thus, the Anderson-Darling estimates can be obtained by solving the following equations simultaneously: $\partial A(\lambda,\alpha)/\partial \lambda =0$ and $\partial A(\lambda,\alpha)/\partial \alpha =0$, where
	\begin{eqnarray}
	\dfrac{\partial A(\lambda,\alpha)}{\partial \lambda} =-\frac{1}{n}\sum_{i=1}^{n}(2i-1)\dfrac{1}{[1-e^{-2\lambda x_{(i)}}]^\alpha}2 x_{(i)}\alpha[1-e^{-2\lambda x_{(i)}}]^{\alpha-1}e^{-2\lambda x_{(i)}} \nonumber \\ +\dfrac{1}{[1-e^{-2\lambda x_{(n+1-i)}}]^\alpha}2 x_{(i)}\alpha[1-e^{-2\lambda x_{(n+1-i)}}]^{\alpha-1}e^{-2\lambda x_{(n+1-i)}}
	\end{eqnarray}
	and
	\begin{equation}
	\dfrac{\partial A(\lambda,\alpha)}{\partial \alpha}=-\frac{1}{n}\sum_{i=1}^{n}(2i-1)2 x_{(i)} \log[1-e^{-2\lambda x_{(i)}}]+ 2 x_{(n+1-i)} \log[1-e^{-2\lambda x_{(n+1-i)}}]
	\end{equation}
\vskip.3cm {\subsection*{2.2.\hskip.3cm Topp Leone $q$ Exponential Distribution}}
\vskip.3cm {\subsection*{2.2.1\hskip.3cm Method of Least Squares Estimation}}
As in the previous case here we can estimate $\lambda$, $\alpha$, and $q$ by minimizing the least square function, with respect to $\lambda$, $\alpha$ and $q$.\\ Therefore LSEs can be obtained by solving the following equations simultaneously:\\
$\dfrac{\partial LS(\lambda,\alpha,q)}{\partial \lambda}=0$ , $\dfrac{\partial LS(\lambda,\alpha,q)}{\partial \alpha}=0,$  $\dfrac{\partial LS(\lambda,\alpha,q)}{\partial q}=0,$ where
\begin{equation}
\dfrac{\partial LS(\lambda,\alpha,q)}{\partial \lambda}= 2\sum_{i=1}^{n}\eta x_{(i)}\alpha (2-q)\psi^{\frac{1}{(1-q)}}[\psi+1]^{\frac{2-q}{1-q}}
\end{equation}
\begin{equation}
\dfrac{\partial LS(\lambda,\alpha,q)}{\partial \alpha} = 2\sum_{i=1}^{n}\eta \ln{[1-\psi^{\frac{2-q}{1-q}}]}\big[1-\psi^{\frac{2-q}{1-q}}\big]^\alpha
\end{equation}
\begin{equation}
\dfrac{\partial LS(\lambda,\alpha,q)}{\partial q}=2\sum_{i=1}^{n}\eta \exp{\big[\ln{[\psi^{\frac{2-q}{1-q}}]}\big]}\alpha [\psi +
1]^{\alpha -1} \big[\frac{\lambda x_{(i)}(2-q)}{(1-q)\psi}+\frac{\ln{[\psi]}}{(1-q)^2}\big]
\end{equation}
where \ \ \ $\psi = 1-\lambda x_{(i)}(1-q)$\ \ \  and $\eta=(1-(1-\lambda x_{(i)}(1-q))^{\frac{2-q}{1-q}})^\alpha -\frac{i}{n+1}$

	\vskip.3cm {\subsection*{2.2.2.\hskip.3cm Method of Weighted Least Squares estimation}}
	Here we can estimate $\lambda$, $\alpha$, and $q$ by minimizing the weighted least square function, with respect to $\lambda$, $\alpha$ and $q$. Therefore LSWEs can be obtained by solving the following equations simultaneously:\\
	$\dfrac{\partial LSW(\lambda,\alpha,q)}{\partial \lambda}=0$ , $\dfrac{\partial LSW(\lambda,\alpha,q)}{\partial \alpha}=0,$  $\dfrac{\partial LSW(\lambda,\alpha,q)}{\partial q}=0,$ where
	\begin{equation}
	\dfrac{\partial LSW(\lambda,\alpha,q)}{\partial \lambda}= 2\sum_{i=0}^{n}\dfrac{(n+1)^2(n+2)}{i(n-i+1)}\eta x_{(i)} \alpha (2-q) \psi^{\frac{1}{1-q}} [1-(1-\psi)^{\frac{1}{1-q}}]^{(\alpha - 1)}
	\end{equation}

	\begin{equation}
	\dfrac{\partial LSW(\lambda,\alpha,q)}{\partial \alpha}= 2\sum_{i=0}^{n}\dfrac{(n+1)^2(n+2)}{i(n-i+1)}\eta x_{(i)} \ln{[1-(1-\psi)^{\frac{2-q}{1-q}}]^\alpha} [1-(1-\psi)^{\frac{2-q}{1-q}}]^\alpha
	\end{equation}
	\begin{eqnarray}
	\dfrac{\partial LSW(\lambda,\alpha,q)}{\partial q}&=& 2\sum_{i=0}^{n}\frac{1}{(q-1)^2(1-\psi)}\big[ e^{\ln{[1-\psi{\frac{2-q}{1-q}}]^\alpha}}\alpha \frac{(n+1)^2(n+2)}{i(n-i+1)}\eta \lambda x_{(i)}\nonumber\\ && \times(2-q) (1-q)
	+\ln(\psi)(\lambda x_{(i)}-\lambda q x_{(i)}-1)(1-(1-\psi)^{\frac{2-q}{1-q}})^{(\alpha-1)}\big]\nonumber\\
	\end{eqnarray}
	where \ \ \ $\psi = 1-\lambda x_{(i)}(1-q)$\ \ \  and $\eta=(1-(1-\lambda x_{(i)}(1-q))^{\frac{2-q}{1-q}})^\alpha -\frac{i}{n+1}$
	
	\vskip.3cm {\subsection*{2.2.3\hskip.3cmMethod of Cramer-von Mises Distance Estimation}}
	In method of Cramer-von Mises Distance Estimation we can estimate $\lambda $, $\alpha$, and $q$ by minimizing the Cramer-von Mises Distance function, with respect to $\lambda$, $\alpha$ and $q$ .\\ Therefore CEs can be obtained by solving the following equations simultaneously:\\
	$\dfrac{\partial C(\lambda,\alpha,q)}{\partial \lambda}=0$ , $\dfrac{\partial C(\lambda,\alpha,q)}{\partial \alpha}=0,$  $\dfrac{\partial C(\lambda,\alpha,q)}{\partial q}=0$ where
\begin{eqnarray}
	\dfrac{\partial C(\lambda,\alpha,q)}{\partial \lambda} &=&  -2 \sum_{i=0}^{n} x_{(i)}\alpha [q\psi^{\frac{1}{1-q}}][(1-\psi)^{(2\alpha-1)}-2(1-\psi^{\frac{1-q}{2-q}}) \nonumber\\ && -\dfrac{q(2i-1)(1-\psi^{\frac{2-q}{1-q}})+2(2i-1)(1-\psi^{\frac{2-q}{1-q}})}{2n}]
	\end{eqnarray}

\begin{equation}
\dfrac{\partial C(\lambda,\alpha,q)}{\partial \alpha}= 2 \sum_{i=0}^{n} \ln{(1-\psi ^{\frac{2-q}{1-q}})} \left[ (1-\psi^{\frac{2-q}{1-q}2 \alpha})  -\dfrac{(2i-1)(1-\psi^{\frac{2-q}{1-q}})^\alpha}{n} \right]
\end{equation}
\begin{eqnarray}
\dfrac{\partial C(\lambda,\alpha,q)}{\partial q} &=& 2 \sum_{i=0}^{n} \big(1-\psi^{\frac{2-q}{1-q}}-\frac{2i-1}{2n} \big)\alpha(1-\psi^{\frac{2-q}{1-q}})^{\alpha-1} e^{\frac{2-q}{1-q}\ln{\psi}}\nonumber\\ &&(\frac{1}{(1-q)^2}\ln{\psi}+\dfrac{\lambda x_{(i)}}{\psi}\dfrac{2-q}{1-q})
	\end{eqnarray}
where \ \ \ $\psi = 1-\lambda x_{(i)}(1-q)$ \ \ \  and  $\eta=(1-(1-\lambda x_{(i)}(1-q))^{\frac{2-q}{1-q}})^{\alpha} -\frac{i}{n+1}$

	\vskip.3cm {\subsection*{2.1.4\hskip.3cm Method of Anderson-Darling Estimation}
		Here we can estimate $\lambda$, $\alpha$, and $q$ of TL$q$E distribution by minimising the Anderson Darling function, with respect to $\lambda$, $\alpha$ and $q$. Therefore Anderson-Darling estimates can be obtained by solving the following equations simultaneously: \\
		$\dfrac{\partial A(\lambda,\alpha,q)}{\partial \lambda}=0$ , $\dfrac{\partial A(\lambda,\alpha,q)}{\partial \alpha}=0,$  $\dfrac{\partial A(\lambda,\alpha,q)}{\partial q}=0$  where \\
	\begin{eqnarray}
		 \dfrac{\partial A(\lambda,\alpha,q)}{\partial \lambda}&=&\dfrac{-1}{n}\sum_{i=1}^{n}(2i-1)
[\dfrac{1}{\xi_{(i)}^{\alpha}}\alpha\xi_{(i)}^{\alpha-1}\dfrac{2(2-q)}{1-q}\psi^{\frac{2(2-q)}{1-q}-1}qx_{(i)}
\nonumber\\&&  -\dfrac{1}{1-\xi_{(n+1-i)}^\alpha}\alpha\xi_{(n+1-i)}^{\alpha-1}\dfrac{2(2-q)}{1-q}\psi^{\frac{2(2-q)}{1-q}-1}qx_{(n+1-i)}]
			\end{eqnarray}

		\begin{equation}
		\dfrac{\partial A(\lambda,\alpha,q)}{\partial \alpha}=\dfrac{-1}{n}\sum_{i=1}^{n}(2i-1)[\log{\xi_{(i)}}+\log{\xi_{(n+1-i)}}]
		\end{equation}

\begin{eqnarray}
		\dfrac{\partial A(\lambda,\alpha,q)}{\partial q}&=&\dfrac{-1}{n}\sum_{i=1}^{n}2(2i-1)\{\dfrac{\alpha \frac{2(2-q)}{1-q} e^{log(\psi_i)}\psi_i \log(\psi_i)+(2-q)(1-q)\lambda x_{(i)}}{\xi_{(i)}(q-1)(1-q)\psi_i} \nonumber \\  &&
		-\dfrac{\alpha \frac{2(2-q)}{1-q} e^{log(\psi_{n+1-i})}\psi_{n+1-i} \log(\psi_{n+1-i})+(2-q)(1-q)\lambda x_{(n+1-i)})}{\xi_{(n+1-i)(1-q)(q-1)\psi_{n+1-i}}}\}\nonumber \\
		\end{eqnarray}
		where $\psi_i = 1-\lambda x_{(i)}(1-q)$ ,\ \ \ $\psi_{n+1-i} = 1-\lambda x_{(i)}(1-q)$ \ \ \ $\xi_{(i)}=[1-[1-\lambda x_{(i)}(1-q)]^{\frac{2(2-q)}{1-q}}]$, \\ $\xi_{(n+1-i)}=[1-[1-\lambda x_{(n+1-i)}(1-q)]^{\frac{2(2-q)}{1-q}}]$
		\paragraph{}
		These are the different estimates of Topp Leone Exponential(TLE) distribution and Topp Leone q Exponential(TLqE) distribution.

%

\end{document}